\documentclass[12pt]{amsart}
\usepackage[margin=1in]{geometry}
\usepackage{amsmath,amssymb,etoolbox,mathrsfs,booktabs,float,enumerate,calc,graphicx}
\usepackage[ruled]{algorithm2e}

\theoremstyle{plain}
\newtheorem{theorem}{Theorem}

\newtheorem*{proposition*}{Proposition}

\newtheorem{lemmma*}{Lemma}

\newtheorem{conjecture}[theorem]{Conjecture}

\theoremstyle{definition}
\newtheorem{definition}[theorem]{Definition}
\newtheorem*{definition*}{Definition}

\newtheorem{remark*}{Remark}

\begin{document}
\title[Norine's Conjecture holds for $n=7$]{Proving Norine's Conjecture holds for $n=7$ via SAT solvers}
\date{\today}
\author[K.\ Frankston]{Keith Frankston}
\author[D.\ Scheinerman]{Danny Scheinerman}
\address{Center for Communications Research, Princeton}
\email{k.frankston@idaccr.org,~~~daniel.scheinerman@gmail.com}
\maketitle
\thispagestyle{empty}

\begin{abstract}
We say a red/blue edge-coloring of the $n$-dimensional cube graph, $Q_n$, is 
\emph{antipodal} if all pairs of antipodal edges have different colors. Norine 
conjectured \cite{Garden} that in such a coloring there must exist a pair of 
antipodal vertices connected by a monochromatic path. Previous work has proven 
this conjecture for $n\le 6$. Using SAT solvers we verify that the conjecture 
holds for $n = 7$.
\end{abstract}

\section{Introduction}
Let $Q_n$ denote the $n$-dimensional cube graph also known as the Hamming cube. 
As is typical, we identify the vertices with binary strings of length $n$. If 
$u$ is a vertex, let $\bar u$ be the vertex with complemented bitstring. We 
say that $u$ and $\bar u$ are antipodal. A pair of edges are antipodal if the set of endpoints of the first are antipodes of the second's. This paper concerns a particular class of edge colorings called antipodal.

\begin{definition}
An \emph{antipodal edge-coloring} of $Q_n$ is a edge $2$-coloring of $Q_n$ such that all pairs of antipodal edges have different colors.
\end{definition}

In this paper we refer to the edge colors as red and blue. We are interested in the following conjecture.

\begin{conjecture}[Norine \cite{Garden}]\label{conj:norine}
For $n\ge 2$, any antipodal edge-coloring of $Q_n$ contains antipodal vertices $u$ and $\bar u$ such that $u$ and $\bar u$ are joined by a monochromatic path.
\end{conjecture}

In Conjecture~\ref{conj:norine} the monochromatic path may be longer than $n$ edges. However, the following (mildly) stronger conjecture seems to hold.

\begin{conjecture}[Geodesic Norine]\label{conj:geodesic}
For $n\ge 2$, for any antipodal edge-coloring of $Q_n$ there exist antipodal vertices $u$ and $\bar u$ such that $u$ and $\bar u$ are joined by a monochromatic geodesic path (path of length $n$).
\end{conjecture}

Early work showed that Conjecture~\ref{conj:norine} holds for $n\le 5$ \cite{FederSubi}. Later, Conjecture~\ref{conj:geodesic} was proven for $n\le 6$ separately in \cite{WestWise} and \cite{CASSAT}. In \cite{CASSAT} the authors use the computer algebra system Sage to make a SAT instance that encodes the space of 
counterexamples for a given $n$. When a solver returns ``UNSAT'' this then 
indicates that there are no counterexamples and therefore the conjecture holds for this value of $n$. Using their 
methodology, they confirm Conjecture~\ref{conj:geodesic} for $n\le 5$ in less 
than 20 seconds of compute time. To verify the case $n=6$ required 1 hour and 35 
minutes for the SAT solver Glucose 3.0. We were unaware of this work when we began 
this project and proceeded in a somewhat similar fashion. However, our encoding 
of possible counterexamples seems significantly more efficient. Using the same 
SAT solver and comparable hardware we are able to establish the case $n=6$ in a 
fraction of a second of compute time. Furthermore, we are able to establish the 
new result that Conjecture~\ref{conj:geodesic} holds for $n=7$ in about three 
minutes of compute time.

We describe our encoding of this problem at a SAT instance in Section~\ref{sec:main} and we give thoughts on how one might approach the case $n=8$ in Section~\ref{sec:eight}.


\section{Main}\label{sec:main}

Our main idea is to encode, for a given $n$, a SAT instance describing a 
counterexample to Conjecture~\ref{conj:geodesic}. Indeed, this project began when the 
second author was skeptical of this conjecture and was seeking a counterexample. 
The encoding is general in the sense that if the system of constraints is 
unsatisfiable then this proves that Norine's conjecture holds for the given $n$. 
As is typical we will describe the problem in conjunctive normal form (CNF). 

We associate to each edge $uv$ of $Q_n$ a binary variable $X_{uv}$ indicating its color (1 is red, 0 is blue). We will have three types of constraints:
\begin{enumerate}
    \item Antipodal edge conditions: Antipodal edges have opposite color.
    \item Geodesic Norine constraints: No pair of antipodal vertices, $u$ and $\bar u$, are connected by a monochromatic geodesic path.
    \item Symmetry breaking.
\end{enumerate}

The antipodal edge constraints are straightforward. There are $n 2^{n-1}$ edges 
in $Q_n$ and consequently $n2^{n-2}$ pairs of antipodal edges. For each pair of 
antipodal edges, $uv$ and $\bar u \bar v$, we have the constraint $X_{uv}+X_{\bar u \bar v}=1$. This can be expressed in CNF as

\begin{equation}\label{eqn:antipodal}
\left(X_{uv} \vee X_{\bar u \bar v}\right) \land \left(\overline{X_{uv}} \vee \overline{X_{\bar u \bar v}}\right).
\end{equation}

Next, we describe the geodesic Norine constraints mentioned above. For each pair of antipodal vertices $u$ and $\bar u$ it suffices to require all connecting geodesic paths contain a red edge. If so then no geodesic path can be entirely red as its complement contains a red edge and thus the original path contains a 
blue edge. So initially we had constraints that for all vertices $u$ and all 
geodesic paths
\begin{equation*}
u:=u_0\sim u_1 \sim \ldots \sim u_n:=\bar u
\end{equation*}
that $\sum_{i=0}^{n-1} X_{u_i u_{i+1}} \ge 1$. This can be improved by noting 
that it suffices to require that the one shorter path $u_0\sim u_1 \sim \ldots 
u_{n-1}$ contains a red edge. This follows as if we let $v=u_{n-1}$ then 
consider the $n+1$ long path
\begin{equation}\label{eqn:bigpath}
 \bar{v}\sim  (u \sim u_1 \sim u_2 \sim \ldots u_{n-2} \sim v)\sim\bar{u}.   
\end{equation}
As $\bar v u$ and $v \bar u$ are antipodal, exactly one of these edges is red. Consequently if the parenthesized path in Equation~\eqref{eqn:bigpath} is entirely blue then there must exist a monochromatic geodesic.

So for each vertex $u$ and each geodesic path $u_0\sim u_1 \sim \ldots \sim u_n$ we will have the constraint $\sum_{i=0}^{n-2} X_{u_i u_{i+1}} \ge 1$. Note that each of these constraints are distinct as the first $n-1$ edges are sufficient to determine a geodesic. As part of our CNF expression we add the clause
\begin{equation}\label{eqn:norine}
    \bigvee_{i=0}^{n-2} X_{u_i u_{i+1}}
\end{equation}
for each geodesic. We have $2^{n-1}$ pairs of antipodal vertices and for each there are $n!$ geodesic paths. So this step adds $2^{n-1} n!$ clauses to our CNF.

Finally, we describe a simple method we used to break the symmetry for this 
problem. These additional clauses are unnecessary for a satisfying arrangement to be a counterexample, but substantially speed up 
the computation. Notice that there must exist a vertex which is incident to 
both red and blue edges. Without loss of generality, it is incident to at least 
as many red edges as blue. So distinguish a vertex $u$ with neighbors $v_1, 
v_2,\ldots v_n$. Let $k = \lceil n/2 \rceil$. Then set $X_{u v_i}=1$ for 
$i=1,\ldots,k$ and $X_{u v_{k+1}}=0$.

Combining these constraints for $n=6$ we generated a SAT instance with $192$ 
variables and $23236$ clauses. For $n=7$ we generated a SAT instance with $448$ 
variables and $323013$ clauses. In addition to Glucose 3.0 we used the solvers 
Kissat, Cryptominisat, Minisat, Parkisat and Clasp. Each of these ran in less 
than a second on the first SAT instance verifying that 
Conjecture~\ref{conj:geodesic} holds for $n=6$. Our new result that this 
conjecture holds for $n=7$ was verified by Glucose 3.0 in 3 minutes and 24 
seconds. Each of the other solvers confirmed this in under six minutes. 

For $n=8$ we generated a SAT instance with $1024$ variables and $5161989$ 
clauses. Early experiments indicated that this would require significant compute 
to solve. The case $n=8$ remains open. We give some thoughts on how one might approach 
this case in Section~\ref{sec:eight}.



\section{Towards $n=8$}\label{sec:eight}

In an attempt to prove Conjecture~\ref{conj:norine} for $n=8$ we looked into more creative ways to break the symmetry. We employed the following theorem \cite{FederSubi}

\begin{theorem}[Feder \& Subi]\label{thm:federsubi}
A $2$-edge-coloring of $Q_n$ is said to by simple if there is no square $xyzt$ 
such that $xy$ and $zt$ are colored red and $yz$ and $tx$ are colored blue. For 
$n\ge 2$, for any simple coloring of $Q_n$ there exists antipodal vertices 
joined by a monochromatic path.
\end{theorem}

As a consequence one can assume without loss of generality that any 
counterexample to Conjecture~\ref{conj:norine} contains a square $xyzt$ such 
that $xy$ and $zt$ are colored red and $yz$ and $tx$ are colored blue. Without 
loss of generality $x,y,z$, and $t$ have bitstrings $00000000$, $01000000$, 
$11000000$ and $10000000$ respectively. These four vertices have an additional 
$4\times 6=24$ edges incident to them. Naively, these edges can be colored in 
$2^{24}$ ways. However, there are three symmetries we apply. Firstly, we can 
apply any of $6!$ permutations to the 3rd through 8th coordinates. Secondly, we 
can swap red versus blue, i.e. complement the 3rd through 8th coordinates. 
Finally, the fixed red and blue edges of the square are preserved under the symmetries
\begin{equation*}
xyzt \to yxtz,~~~~~xyzt \to tzyx,~~~~~xyzt \to ztxy.
\end{equation*}
Under these symmetries the $2^{24}$ possible colorings are reduced to 7218 orbits\footnote{We enumerated these directly via Julia code, but also confirmed this count via applying Burnside's lemma}. This means we can reduce resolving Conjecture~\ref{conj:norine} to 7218 subproblems where we have 28 edges colored, i.e. 28 variables get an assigned value. We hoped these 
subproblems would be manageable (maybe requiring a few minutes of compute time 
on average). We chose 100 of these 7218 subproblems at random and ran them single threaded using the SAT solver Kissat on a 3.6 GHz cpu. On average they ran in 2.9 days with the longest taking 4.7 days. Thus we make the (very rough) estimate that this approach would require 
\begin{equation*}
7218 \times 2.9~\mbox{cpu days} = 57.3~\mbox{cpu years}.
\end{equation*}
We are optimistic that with advances in speed of SAT solvers and hopefully a few 
additional good ideas this computation will become more readily doable.



\end{document}